\input amstex
\documentstyle{amsppt}
\input symm3.def
\onefiletrue
\finaltrue

\topmatter
\title The Point Spectrum of the Dirac Operator on Noncompact Symmetric
Spaces\endtitle
\author S. Goette, U. Semmelmann\endauthor
\rightheadtext{the Dirac Operator on Noncompact Symmetric Spaces}
\address  Universit\'e de Paris-Sud, D\'epartement de Math\'ematique,
Laboratoire de Topologie et Dynamique, URA D1169 du CNRS, B\^atiment 425,
91405 Orsay Cedex, France\endaddress
\email Sebastian.Goette\@math.u-psud.fr\endemail
\address Mathematisches Institut, Universit\"at M\"unchen,
Theresienstr.~39, D-80333 M\"unchen, Germany\endaddress
\email semmelma\@rz.mathematik.uni-muenchen.de\endemail
\thanks Both authors were supported by a research fellowship of the~DFG
\endthanks
\abstract
In this note, we consider the Dirac operator~$D$
on a Riemannian symmetric space~$M$ of  noncompact type.
Using representation theory,
we show that~$D$
has point spectrum iff the $\Adach$-genus of its compact dual
does not vanish.
In this case, if~$M$ is irreducible
then $M=\U(p,q)/\U(p)\times\U(q)$ with~$p+q$ odd,
and~$\Spec_p(D)=\{0\}$.
\endabstract
\endtopmatter

\document
\head 0. Introduction\endhead
We investigate the existence of point spectrum of the Dirac operator~$D$
acting on spinors over a Riemannian symmetric space~$M=G/K$ of noncompact type.
Following Seifarth's approach in~\Seifarth,
we look at those discrete series representations of~$G$
that appear in~$L^2(\calS)$,
where~$\calS$ is the spinor bundle over~$M$.
We find that the existence of point spectrum of~$D$
is equivalent to a regularity condition for the half sum~$\rho_\frk$
of positive roots of~$K$,
which in turn is equivalent to the nonvanishing of the $\Adach$-genus
of the compact dual~$M'$ of~$M$.
Using 
the classification of compact symmetric spaces,
we finally determine all noncompact symmetric spaces
on which the Dirac operator~$D$ has point spectrum.
We summarize our results:

\Theorem\MainTheorem.
Let~$M$ be a Riemannian symmetric space of noncompact type,
and let~$D$ be the Dirac operator acting on spinors over~$M$.
Then the following statements are equivalent:
\roster
\item the point spectrum of~$D$ is nonempty;
\item the point spectrum of~$D$ is precisely~$\Spec_p(D)=\{0\}$;
moreover, as a $G$-module,
$\ker(D)$ is irreducible and isomorphic to the
discrete series representation with Harish-Chandra parameter~$\rho_\frk$;
\item the $\Adach$-genus of the compact dual of~$M$ is nonzero;
\item each irreducible factor of~$M$
is isometric to~$\U(p,q)/\U(p)\times\U(q)$,
with~$p+q$ odd.
\endroster

Our present note is motivated by the work of several authors on the spectra
of Dirac operators on noncompact Riemannian symmetric spaces. 
Using the Plancherel theorem,
Bunke computed the whole spectrum of the untwisted Dirac operator~$D$
on the real hyperbolic spaces in~\Bunke\
(note the incorrect statement concerning the eigenvalue~$0$).
Seifarth showed the nonexistence of point spectrum on the real and
quaternion hyperbolic spaces in~\Seifarth\
(the treatment of the complex hyperbolic space is incomplete).
Another computation of the spectrum of~$D$ on~$\R H^n$ by Camporesi and Higuchi
uses polar coordinates and separation of variables (\CH).
Using a similar approach,
Baier proved in~\Baier\ that the Dirac operator on~$\C H^n$ has no
eigenvalue~$\lambda$ with~$\abs\lambda\ge{n-1\over4}$.
%
Let us also mention the results of Galina and Vargas on the eigenvalues
of twisted Dirac operators:
In~\GV,
they compute the spectrum of Dirac operators on~$\R H^n$ and~$\C H^n$,
twisted with a homogeneous vector bundle.
They consider only the case
where the inducing $K$-representation has a  sufficiently nonsingular
highest weight.

The rest of this paper is organized as follows:
In \forward\ASKapitel,
we recall the relation between point spectrum of
homogeneous selfadjoint elliptic operators on~$M=G/K$
and discrete series representations of~$G$.
In \forward\KTypKapitel,
we show that the existence of point spectrum of~$D$ on~$M$
is equivalent to the nonvanishing of the $\Adach$-genus
on the compact dual of~$M$.
Finally, in \forward\BeispielKapitel,
we classify the compact symmetric spaces~$M'$ with~$\Adach(M')[M']\ne0$.

This work was written while the second named author enjoyed the hospitality
and support of the~IHES (Bures-sur-Yvette).
The first named author would like to thank
the Universit\'e de Paris-Sud (Orsay) for its hospitality.
We are grateful to C.~B\"ar, J.-M.~Bismut and W.~M\"uller
for helpful discussions.
We wish to thank M.~Olbrich for carefully reading the manuscript,
pointing out a few inaccuracies,
and suggesting an alternative proof of \AdachTheorem.

\head 1. The Point Spectrum and the Discrete Series\endhead
\Kapitel\ASKapitel=1

Let~$M = G/K$  be  a Riemannian symmetric space of noncompact type.
Here, $G$  is a noncompact connected semisimple Lie group,
and~$K$ is a maximal compact subgroup.
We fix a $G$-invariant metric on~$M$.
Then~$M$ is a Hadamard manifold, 
i.e.\ the Riemannian exponential map~$\exp\colon T_pM\to M$
is a diffeomorphism at each point~$p$ of~$M$.
In particular, $M$ is contractible,
and thus possesses a unique spin structure.

Let~$\frg=\frk\oplus\frp$ be the Cartan
decomposition of~$\frg$, where~$\frk$ is
the Lie algebra of~$K$.
A homogeneous spin structure can be described by a lift~$\tilde \alpha $
of the adjoint representation~$\alpha\colon K \to \SO(\frp)$
to~$\Spin(\frp)$.
We can assume the existence of such a lift
(if necessary, we replace~$G$ and~$K$ by suitable double covers).
The (complex) spin representation~$(\rho, S)$ of~$\Spin(\frp)$ gives rise
to a $K$-representation~$(\sigma, \, S)$,
with~$\sigma:=\rho \circ {\tilde \alpha} $.
The spinor bundle is then isomorphic to the homogeneous vector bundle
	$$\calS:= G \times_\sigma S\Formel\SpinorBundleDefinition$$
induced by~$\sigma$.
The Levi-Civita connection on~$M$ induces a connection on~$\calS$.

Let $\Gamma_c(\calS)$ be the space of compactly supported smooth
sections of $\calS$, and let~$L^2(\calS)$ be its Hilbert 
space completion.
The Dirac operator acts on~$\Gamma_c(\calS)$ as the composition
of covariant derivative and Clifford multiplication.
Since $M = G/K$  is a complete manifold, the Dirac operator is
essentially selfadjoint (cf. \Wolf).
Hence, its minimal
and maximal closed extension coincide.
Let~$D$ be the unique selfadjoint extension to a closed operator.
It commutes with the natural action of~$G$ on~$L^2(\calS)$.
More generally,
one can consider Dirac operators on~$L^2(\calS\tensor\calW)$,
where~$\calW$ is a homogeneous Hermitian vector bundle over~$M$
which is equipped with an equivariant unitary connection.

Because~$D$ is selfadjoint,
its spectrum 
consists only of point spectrum and continuous spectrum.
Moreover, it
is completely contained in~$\R$.
The {\em point spectrum\/}~$\Spec_p(D)$, i.e. the set of eigenvalues, is
defined as
	$$\Spec_p(D)
	:=\bigl\{\,\lambda\in\C\bigm|\ker(D-\lambda)\neq\{0\}\,\bigr\}\;.$$
If~$\lambda$ is an eigenvalue of~$D$,
the dimension of the eigenspace~$\ker(D - \lambda)$
is called the {\em multiplicity\/} of~$\lambda$.

Clearly, $G$ acts on the eigenspaces of~$D$.
It turns out that the eigenspaces are direct sums
of irreducible $G$-representations belonging to the discrete series:

\Definition\DiscreteSeriesDefinition
An irreducible representation~$(\pi,H)$ of~$G$ is called a 
{\em discrete series representation\/}
iff the matrix coefficients~$g\mapsto\<\pi(g) v, w\>$
for all~$v$, $w\in H$
are square integrable on~$G$
with respect to the Haar measure.
Let~$\Gdach_d$ be the set of equivalence classes
of discrete series representations.
\enddefinition

The main tool of our investigation of~$\Spec_p(D)$
is the following

\Theorem\PointSpectrumTheorem{(cf.~\AS, \CM)}.
Let~$D$ be a homogeneous selfadjoint elliptic differential operator
on~$\calE:=G\times_\eps E$ for some $K$-representation~$(\eps,E)$.
Then the direct sum of all eigenspaces of~$D$
is isomorphic to
	$$\bigoplus_{\pi\in\Gdach_d}\pi\tensor\Hom_K\(\pi|_K,\eps\)\;.$$

In particular, a discrete series representation~$\pi$ of~$G$
is isomorphic to a subrepresentation of~$L^2(\calE)$
iff~$\pi|_K$ has an irreducible $K$-subrepresentation
in common with~$\eps$.
In this case, we say that~$\pi\in\Gdach_d$ contributes to~$\Spec_p(D)$.

\demo{Proof\/ \rm of \PointSpectrumTheorem}
Since~$D$ is a $G$-invariant elliptic differential operator,
we can apply a theorem of Connes and Moscovici
(\CM, Theorem~6.1). It follows that each eigenspace of~$D$ is
isomorphic to a finite sum of discrete series representations of~$G$.

On the other hand,
by the Plancherel Theorem and Frobenius reciprocity
(cf.~\AS, chapter~2),
we have
	$$\Hom_G\(\pi,L^2(\calE)\)\cong\Hom_K\(\pi|_K,\eps\)$$
for each discrete series representation~$\pi$.
Moreover,
$D$ is $G$-invariant and~$\Hom_K\(\pi|_K,\eps\)$
is finite dimensional by results of Harish-Chandra
(Theorem~8.1 in~\Knapp).
Hence,
it is easy to check that
	$$\pi\tensor\Hom_G\(\pi,L^2(\calE)\)\subset L^2(\calE)$$
decomposes as a finite sum of $D$-eigenspaces.
\qed\enddemo

Each eigenvalue has infinite multiplicity,
since all nontrivial unitary representations of a noncompact connected
semisimple Lie group are infinite dimensional.
Moreover, if $D$  has nonempty point spectrum
then~$G$ has discrete series representations.
Due to a theorem of Harish-Chandra
(cf.~\Knapp, Theorem~12.20, \AS, Proposition~6.11),
this is the case iff~$\rk(G) =\rk(K)$.
Hence, we have the following

\Remark\EqualRankRemark
On a noncompact symmetric space~$G/K$ with~$\rk(G)>\rk(K)$,
the point spectrum of the Dirac operator~$D$ is empty.
\endremark

\head 2. Minimal $K$-Types and Point Spectrum\endhead
\Kapitel\KTypKapitel=2
In this chapter,
we recall a few facts from the theory
of discrete series representations of~$G$.
We will show that at most one irreducible subrepresentation of~$\sigma$
can occur as a $K$-type of a discrete series representation of~$G$.
This happens iff the half sum~$\rho_\frk$ of positive roots of~$K$
is $\frg$-regular.
Using this,
we prove the equivalence of statements~(1) --~(3) of \MainTheorem.
We remark that our arguments in this chapter are also valid
for nonirreducible symmetric spaces.

\subhead a) Discrete Series Representations and their $K$-Types\endsubhead
By \EqualRankRemark, we may and will assume from now on
that~$\rk(G)=\rk(K)$.
Then we fix a common maximal torus~$H\subset K\subset G$
with Lie algebra~$\frh$ and weight lattice
	$$\Gamma:=\bigl\{\,\gamma\in i\frh^*\bigm|\gamma(X)\in2\pi i\Z
		\hbox{ for all~$X\in\frh$ with }e^X=e\,\bigr\}\;.$$
Let~$\Delta_\frg=\Delta_\frk\cup\Delta_\frp$
be the root system of~$G$ with respect to~$\frh$,
decomposed into the root system of~$K$
and the set of noncompact roots.
Choose systems of positive roots~$\Delta_\frg^+\supset\Delta_\frk^+$,
and let~$P_\frg\subset P_\frk\subset i\frh^*$ be the Weyl chambers
associated to~$\Delta_\frg^+$ and~$\Delta_\frk^+$.
Let~$W_\frg$, $W_\frk$ be the Weyl groups of~$G$ and~$K$,
and set
	$$W':=\{\,w\in W_\frg\mid w(P_\frg)\subset P_\frk\,\}\;.$$
Let~$\rho_\frg$ and~$\rho_\frk$ be the half sums of positive roots
of~$G$ and~$K$.

We fix an $\Ad_K^*$-invariant scalar product on~$\frk^*$.
We call a weight~$\lambda\in i\frh^*$ {\em $\frg$-regular\/}
if~$\<\lambda,\gamma\>\ne 0$ for all~$\gamma\in\Delta_\frg$,
and {\em $\frg$-singular\/} otherwise.
If~$\<\lambda,\gamma\>\ne 0$ holds only for~$\gamma\in\Delta_\frk$,
then~$\lambda$ is called {\em $\frk$-regular.\/}
Clearly, $\frg$-regularity implies $\frk$-regularity.
An element~$\kappa\in i\frh^*$ is called {\em $\frk$-algebraically
integral\/} iff
	$$2{\<\alpha,\kappa\>\over\<\alpha,\alpha\>}\in\Z$$
for all~$\alpha\in\Delta_\frk^+$.
Note that all weights of~$K$, i.e.\ all elements of~$\Gamma$,
are automatically $\frk$-algebraically integral.
Note also that~$\rho_\frk$ and~$\rho_\frg$ are $\frk$-algebraically integral
(for~$\rho_\frk$ this is well known,
for~$\rho_\frg$ it follows because~$\rho_\frg$ is
$\frg$-algebraically integral).
Furthermore, $\rho_\frk$ uniquely minimizes~$\abs\kappa$ among all
$\frk$-algebraically integral $\frk$-regular elements~$\kappa$ of~$P_\frk$
(for semisimple~$K$, this is well known,
in the general case it follows because the center of~$K$ is orthogonal
to its semisimple part).

Let us now turn to some facts about the discrete series of~$G$,
in particular about the possible $K$-types.

\Definition\MiniKTypDef
Let~$\pi\in\Gdach_d$ be a discrete series representation of~$G$,
and let~$\phy_\kappa$ be an irreducible representation of~$K$
with highest weight~$\kappa\in\Gamma\cap P_\frk$.
Then~$\kappa$ is called a {\em $K$-type of~$\pi$,\/}
if~$\pi|_K$ contains an irreducible subrepresentation
isomorphic to~$\phy_\kappa$.
The dimension of~$\Hom_K(\pi,\phy_\kappa)$ is called the {\em multiplicity\/}
of~$\kappa$.

If~$\kappa$ minimizes~$\left|\kappa+2\,\rho_\frk\right|$ among all $K$-types,
then~$\kappa$ is called a {\em minimal $K$-type of~$\pi$.}
\enddefinition

Our argumentation is based upon the following 
fundamental result of Harish-Chandra:

\Theorem\KTypeTheorem({\AS, Theorems~8.1 and~8.5,
\Knapp, Theorems 9.20 and 12.21}).
The discrete series representations of~$G$ are parametrized
by~$\lambda\in P_\frk$ with~$\lambda\in w(P_\frg)$ for some~$w\in W'$
such that~$\lambda$ is regular
and~$\lambda - w\rho_\frg\in\Gamma$.
For such a~$\lambda$,
the discrete series representation~$\pi_\lambda$
corresponding to~$\lambda$ has a unique minimal $K$-type
	$$\kappa:=\lambda+w\rho_\frg-2\,\rho_\frk\;,$$
which occurs with multiplicity~$1$.
Finally, each $K$-type~$\kappa'$ of~$\pi_\lambda$ is of the form
	$$\kappa':=\kappa \, + \sum_{\textstyle{\alpha\in\Delta_\frg
			\atop\<w\rho_\frg,\alpha\>>0}}n_\alpha\,\alpha
	\Formel\KTypFormel$$
where the~$n_\alpha$ are nonnegative integers.

In the literature, $\lambda$ is called the {\em Harish-Chandra parameter\/}
for~$\pi_\lambda$,
while~$\kappa$ is called the {\em Blattner parameter.\/}

\subhead b) The Dirac Operator on Spinors\endsubhead
Let~$\calS=G\times_\sigma S$ be the spinor bundle on~$M$
as in~\SpinorBundleDefinition.
Since we assume that~$\rk(G)=\rk(K)$,
the symmetric space~$M$ is even dimensional.
In particular, the spinor representation and the spinor bundle split
into a positive and a negative part:
	$$S=S^+\oplus S^-\;,\qquad\text{and}\qquad
	\calS=\calS^+\oplus\calS^-\;.$$
The $K$-action on~$S^+$ and~$S^-$ is described
by a formula of Parthasarathy:

\Lemma\ParthLemma{(\forward\Parth, Lemma~2.2)}.
For~$w\in W'$,
let~$\sigma^{w\rho_\frg-\rho_\frk}$
be the $K$-representation with highest weight~$w\rho_\frg-\rho_\frk\in P_\frk$.
Then for a suitable orientation of~$M$,
$\sigma$ decomposes as
	$$\sigma=\sigma^+\oplus\sigma^-
		:=\bigoplus_{\textstyle{w\in W'\atop\sign(w)=1}}
			\sigma^{w\rho_\frg-\rho_\frk}
		\oplus\bigoplus_{\textstyle{w\in W'\atop\sign(w)=-1}}
			\sigma^{w\rho_\frg-\rho_\frk}\;.$$

\Remark\IntegralityRemark
This implies in particular that~$\rho_\frk - w\rho_\frg\in\Gamma$,
because we have assumed that~$\sigma$ is a representation of~$K$.
Note that~$w\in W_\frg$ may be arbitrary,
because different $W_\frg$-translates differ by linear combinations of roots
of~$G$,
which are clearly in~$\Gamma$.
\endremark

\Remark\PlusMinusRemark
We recall that the operator~$D$ splits as
	$$D^\pm:=D|_{\Gamma(\calS^\pm)}
		\colon\Gamma(\calS^\pm)\to\Gamma(\calS^\mp)\;.$$
If~$E_\mu$ is an eigenspace corresponding to an eigenvalue~$\mu$ of~$D$,
then~$E_\mu$ splits into~$E_\mu^+\oplus E_\mu^-$,
with~$E_\mu^\pm:=E_\mu\cap\Gamma(\calS^\pm)$.
If moreover, $\mu\ne 0$,
then~$D^\pm|_{E_\mu^\pm}\colon E_\mu^\pm\to E_\mu^\mp$ is an isomorphism.
\endremark

We will now establish an algebraic criterion for the existence
of point spectrum for untwisted Dirac operators.

\Theorem\DiracTheorem.
Let~$D$ be the untwisted Dirac operator on~$M=G/K$.
If~$\rho_\frk$ is $\frg$-regular,
then~$\Spec_p(D)=\{0\}$,
and~$\ker(D)$ is isomorphic to the discrete series representation
with Harish-Chandra parameter~$\rho_\frk$.
If~$\rho_\frk$ is $\frg$-singular,
then there is no point spectrum.

\Rem
By~\forward\AS, Theorem~9.3,
we already know that~$\ker(D)\ne 0$ iff~$\rho_\frk$
is $\frg$-regular.
It would thus be enough to check that~$\Spec_p(D)$ contains
no nonzero eigenvalues.
\endremark

\Proof
First of all,
if~$\rho_\frk$ is $\frg$-regular,
then there exists a discrete series representation
with Harish-Chandra parameter~$\rho_\frk$
because of \KTypeTheorem\ and~\IntegralityRemark.
Let~$w\in W'$ such that~$\rho_\frk\in w(P_\frg)$.
The minimal $K$-type of~$\pi_{\rho_\frk}$ is~$w\rho_\frg-\rho_\frk$,
which is a highest weight of~$\sigma$ by~\ParthLemma.
Hence by \PointSpectrumTheorem, $D$ has point spectrum.

On the other hand,
let us assume that~$\pi_\lambda$ is a discrete series representation of~$G$
that contributes to~$\Spec_p(D)$.
We will show that then necessarily~$\lambda=\rho_\frk$.

Let~$w\in W'$ be such that~$\lambda\in w(P_\frg)$.
By \PointSpectrumTheorem\ and~\ParthLemma,
for some~$w_0\in W'$,
the weight~$w_0\rho_\frg-\rho_\frk$ is a $K$-type of~$\pi_\lambda$.
Then by \KTypeTheorem,
there exist nonnegative integers~$n_\alpha$ such that
	$$w_0\rho_\frg+\rho_\frk
	=\lambda+w\rho_\frg
		+\sum_{\textstyle{\alpha\in\Delta^+_\frg
		\atop\<w\rho_\frg,\alpha\>\ge 0}}n_\alpha\,\alpha\;.
	\Formel\DiracKTypFormel$$

We establish a few inequalities:
By construction, $\lambda$ and~$w\rho_\frg$ are both $\frg$-regular
and lie in the same Weyl chamber~$w\(P_\frg\)$ of~$\frg$.
This has two consequences:
First,
the weight~$w\rho_\frg$ uniquely minimizes the distance to~$\lambda$
among all $W_\frg$-translates of~$\rho_\frg$.
Thus
	$$\<\lambda,w_0\rho_\frg\>\le\<\lambda,w\rho_\frg\>\;,
	\Formel\UnglEins$$
with equality iff~$w_0=w$.
Second,
$\<\alpha,w\rho_\frg\>>0$ iff~$\<\alpha,\lambda\>>0$.
This implies
	$$\left\<\lambda,\sum
		n_\alpha\,\alpha\right\>\ge0\;,
	\Formel\UnglZwei$$
since the~$n_\alpha$ have to be nonnegative.
Moreover, we have equality iff all the~$n_\alpha$ are zero.

If~\DiracKTypFormel\ holds for~$\lambda$,
then~$\lambda$ must clearly be $\frk$-algebraically integral,
because the same holds for~$\rho_\frk$, $\rho_\frg$ and all~$\alpha\in\Delta^+_\frg$.
Now the weight~$\rho_\frk$ uniquely minimizes~$\abs\kappa$
among all $\frk$-regular $\frk$-algebraically integral~$\kappa\in P_\frk$.
Because~$\lambda$ is $\frg$-regular, it is also $\frk$-regular,
and we have
	$$\<\rho_\frk,\lambda\>
		\le\abs\lambda\cdot\abs{\rho_\frk}
		\le\abs\lambda^2
	\Formel\UnglDrei$$
with equality iff~$\lambda=\rho_\frk$.

In order to show that~$\lambda=\rho_\frk$,
we multiply~\DiracKTypFormel\ by~$\lambda$ and apply~\UnglEins\ und~\UnglZwei:
	$$\left\<w_0\rho_\frg+\rho_\frk,\lambda\right\>
	=\left\<\lambda+w\rho_\frg+\sum
		n_\alpha\,\alpha,\lambda\right\>
	\qquad\Longrightarrow\qquad\qquad
	\left\<\rho_\frk,\lambda\right\>
	\ge\abs\lambda^2\;.$$
So by~\UnglDrei,
we have equality, which means that~$\lambda=\rho_\frk$, that~$w_0=w$,
and that all~$n_\alpha$ are zero.
Now, $\rho_\frk$ can only be a Harish-Chandra parameter for a discrete series
representation of~$G$ if $\rho_\frk$ is $\frg$-regular.
Thus, there is no point spectrum if~$\rho_\frk$ is $\frg$-singular.

Let us assume that~$\rho_\frk$ is $\frg$-regular.
Then, among the highest weights of~$\sigma$,
only~$w\rho_\frg-\rho_\frk$ can appear as a $K$-type of a discrete
series representation of~$G$.
This implies that the eigenspaces of~$D$ are contained
either in~$L^2(\calS^+)$ or in~$L^2(\calS^-)$.
In particular, $\Spec_p(D)\subset\{0\}$,
because by \PlusMinusRemark, any nonzero eigenvalue~$\mu$
would lead to an eigenspace~$E_\mu=E_\mu^+\oplus E_\mu^-$
with~$E_\mu^+\cong E_\mu^-\not\cong\{0\}$.
Finally,
$w\rho_\frg-\rho_\frk$ is the minimal $K$-type of~$\pi_{\rho_\frk}$.
Hence, it has multiplicity~$1$,
and~$\ker(D)$ is irreducible as a $G$-module.
\qed\enddemo

\Rem
Another way to check that~$D$ vanishes
on~$\pi_{\rho_\frk}\subset L^2(\calS)$
is to express~$D^2$ in terms of the Casimir operator~$\Omega$ of~$G$
(\Parth, Proposition~3.1, \Knapp, Lemma~12.12),
and using the explicit formula for~$\pi_\lambda(\Omega)$
(\Knapp, Lemma~12.28).
\endremark\medskip

We will now reformulate the theorem above in terms of the compact dual
of~$M$.
Therefore, let~$\frg^\C$ be the complexification of~$\frg$.
Recall that there exists a compact, connected, simply connected Lie group~$G'$
with Lie algebra~$\frg':=\frk\oplus i\frp$.
Let~$K'\subset G'$ be its Lie subgroup with Lie algebra~$\frk$,
then~$K'$ is closed,
and~$M':=G'/K'$ is called the {\sl compact dual of~$M$.\/}
Note that~$\frh$ is a common Cartan subalgebra of~$\frg$, $\frg'$ and~$\frk$,
and that~$\frg$ and~$\frg'$ have the same roots, Weyl chambers etc.
with respect to~$\frh$.

With these definitions, we can give an equivalent criterion for the
existence of point spectrum:

\Corollary\AdachTheorem.
Let~$D$ be the untwisted Dirac operator on~$M=G/K$.
Then~$D$ has point spectrum iff the $\Adach$-genus of the compact
dual~$M'=G'/K'$ of~$M$ is nonzero.

\Proof
By \forward\BH,
the $\Adach$-genus of~$M'$ is given by the formula
	$$\Adach(M')[M']
	=\prod_{\alpha\in\Delta_\frg}
		{\<\alpha,\rho_\frk\>\over\<\alpha,\rho_\frg\>}
	\Formel\ADachFormel$$
for a suitable orientation of~$M'$.
In particular, $\Adach(M')[M']\ne 0$ iff~$\rho_\frk$ is $\frg$-regular,
cf.\ Theorem~23.3 in~\BH.
Thus, our claim follows from \DiracTheorem.
\qed\enddemo

\Rem
The following alternative proof motivates the appearance of the $\Adach$-genus:
By \DiracTheorem,
$D$ has point spectrum iff its $L^2$-index is nonzero.
Using Hirzebruch proportionality,
one then concludes that this is the case iff~$\Adach(M')[M']\ne0$
(cf.~\AS, chapter~3 and erratum).
This was suggested by M.~Olbrich.
\endremark

Let us state some consequences of our criterion.

\Corollary\ComplexCorollary.
Let~$D$ be the untwisted Dirac operator on~$M=G/K$.
If~$D$ has point spectrum, then
\roster
\item $M$ is a Hermitian symmetric space,
\item the compact dual~$M'$ of~$M$ carries no spin structure, and
\item the dimension of~$M$ is divisible by~$4$.
\endroster

\Proof
By \BH, Theorem~23.3, $M'$ is Hermitian symmetric if its $\Adach$-genus
is nonzero.
Then~$M$ is also Hermitian symmetric.
Next, $M'$ has positive scalar curvature.
Thus if~$M'$ was spin,
its $\Adach$-genus would vanish by Lichnerowicz' theorem
(\LM, Corollary~8.9).
Finally, the $\Adach$-genus of~$M'$ can be nonzero only if~$\dim M'$
is divisible by~$4$.
Hence, the claims follows from \AdachTheorem.
\qed\enddemo

\Rem
The conditions listed in \ComplexCorollary\ are not sufficient for the
existence of point spectrum:
In the next section,
we will see that for~$M':=\SP(n)/\U(n)$ with~$n\in4\N$,
conditions~(1) --~(3) above are satisfied.
Nevertheless, the $\Adach$-genus of~$M'$ vanishes.
\endremark

\head 3. Compact Symmetric Spaces with Nonvanishing $\Adach$-Genus\endhead
\Kapitel\BeispielKapitel=3
In this section we want to determine the compact Riemannian
symmetric spaces~$M'=G'/K'$
with non-vanishing $\Adach$-genus.
By~\Bott, we may again assume that~$\rk(G') = \rk(K')$.
Because the $\Adach$-genus is multiplicative on products of manifolds,
we restrict our attention to irreducible symmetric spaces.
By \ComplexCorollary,
we only have to investigate compact Hermitian symmetric spaces~$M'$
with~$\dim M'\in4\N$ which are not spin.

The simply connected symmetric spaces that admit a spin
structure are known (cf.~\CG\ or~\HS).
On the other hand there are four families of Hermitian symmetric spaces
and two exceptional ones (cf.~\Hel).
Combining these lists,
we see that the irreducible Hermitian symmetric spaces
which have no spin structure
form the following three families:
$$\alignat 3
	(1)&\qquad&
		\SO(n + 2)&/\SO(2)\times\SO(n)\qquad&
			&\text{for~$n$ odd,}\\
	(2)&&
		\SP(n)&/\U(n)&
			&\text{for~$n$ even,}\quad\text{and}\\
	(3)&&
		\U(p + q)&/\U(p)\times\U(q)&
			&\text{for~$p+q$ odd}
\endalignat$$
(in the last case, $M$ should actually be represented
as a quotient of a finite cover of~$\SU(p,q)$,
rather than of~$\U(p,q)$).
We will show that all manifolds of the families~(1) and~(2)
have vanishing $\Adach$-genus,
while the $\Adach$-genus of~$\U(p + q)/\U(p)\times\U(q)$
is different from zero if~$p+q$ is odd.

The manifolds of family~(1) have dimension~$2n$.
Since we assume that~$n$ is odd, the dimension is not divisible by~$4$,
and the $\Adach$-genus vanishes.
For the two other families, we have to use formula~\ADachFormel\
and to compute the scalar products~$\<\alpha,\rho_{\frk}\>$
for all positive roots~$\alpha\in\Delta_\frg^+$.

Let us investigate family~(2).
The Lie algebra of~$\SP(n)$
is the Lie algebra of skew-Hermitian quaternionic matrices of order~$n$. 
The Lie algebra of~$\U(n)$ is realized
as the sub-algebra of skew-Hermitian complex matrices of order~$n$.
A common Cartan sub-algebra~$\frh$ is the Lie algebra of the matrices
of the form
	$$\lambda=(\lambda_1, \ldots, \lambda_n)
	:=\diag( \lambda_1, \ldots, \lambda_n)\;,$$
with~$\lambda_j\in i\R$.
As a system of positive roots of~$\SP(n)$, we take
	$$\Delta_\frg^+
	=\{\,\lambda_i\pm\lambda_j\mid 1\le i < j\le n\,\}
		\cup\{\,2\lambda_i\mid i=1, \dots, n\,\}\;.$$
The positive roots in $K' = \U(n)$ are
$\Delta_\frk^+:=\{\,\lambda_i - \lambda_j\mid i < j\,\}$.
Hence, $\rho_{\frk}=(n-1, n - 3, \dots , 1-n)$.
Clearly,
the standard scalar product on~$\frh^*\cong\R^n$
extends to an $\Ad_{G'}^*$-invariant scalar product on~$\frg^{\prime*}$.
In particular,
for~$\alpha=(\lambda_1 + \lambda_n)$ we have~$\<\alpha,\rho_{\frk}\>=0$.
Hence, according to formula~\ADachFormel,
the $\Adach$-genus of all manifolds~$\SP(n)/\U(n)$ is zero.

The computation for family~(3) is similar.
As a system of positive roots we take
	$$\Delta_\frg^+
	=\{\lambda_i-\lambda_j\mid 1 \le i < j \le p + q\,\}\;.$$
The positive roots in $K' = \U(p) \times \U(q)$
are
	$$\Delta_\frk^+
	=\{\,\lambda_i-\lambda_j\mid 1 \le i < j \le p\,\}
		\cup\{\,\lambda_i-\lambda_j\mid p+1\le i<j\le p+q\,\}\;.$$
This yields~$\rho_{\frk}=(p-1, p-3, \ldots, 1-p, q-1, q-3, \ldots 1-q)$. 
Again we take as a scalar product for the roots the canonical scalar product
of vectors in~$\R^n$.
Since~$p+q$ is odd,
we can assume~$p$ to be even and~$q$ to be odd.
Hence, all numbers~$p-1$, $p-3$, \dots, $1-p$ are odd and all numbers
$q-1$, $q-3$, \dots, $1-q$ are even. 
From this it follows that the scalar product~$\<\alpha,\rho_{\frk}\>$
for any positive root~$\alpha\in\Delta_\frg^+$ is different from zero. 
Using once again formula~\ADachFormel,
we obtain that the $\Adach$-genus of~$\U(p + q)/\U(p) \times \U(q)$
is nonzero.

In particular, the $\Adach$-genus of the complex projective space~$\C P^{2n}$
does not vanish.
Here a simple computation
gives~$\Adach(\C P^{2n})[\C P^{2n}]=(-4)^{-n}\,\prod^n_{i=1}{2i-1\over2i}
=(-16)^{-n}\bigl({2n\atop n}\bigr)$.
Finally, we have

\Theorem\Classification.
Let $ M' = G'/K' $ be an irreducible Riemannian symmetric space
of compact type.
Then $M'$ has nonvanishing $\Adach$-genus iff~$M'$ is isometric to
	$$\U(p + q)/\U(p) \times \U(q)\;,
		\qquad\text{with $p + q$ odd.}
	\quad\qed$$

Together with \DiracTheorem, \AdachTheorem,
and the multiplicativity of~$\Adach$,
this proves~\MainTheorem.
\qed

\Refs
\widestnumber\key{CM}

\Quelle\AS[AS]\Journal
  M.~F.~Atiyah, W.~Schmid:
  A Geometric Construction of the Discrete Series
  for Semisimple Lie Groups,
  Inv.\ math.~42 (1977), 1--62
  \moreref\paper Erratum
  \jour Inv.\ math.\vol 54 \yr 1979\pages 189--192

\Quelle\Baier[Ba]\Preprint
P.~D.~Baier:
  \"Uber den Diracoperator auf Mannigfaltigkeiten mit Zylinderenden,
  Diplomarbeit, Universit\"at Freiburg \yr 1997

\Quelle\Baer[B\"a]\Journal
  C.\ B\"ar:
  The Dirac fundamental tone of the hyperbolic space,
  Geom. Dedicata~41 (1992), 103--107

\Quelle\BGV[BGV]\Buch
  M.~Berline, E.~Getzler, N.~Verne:
  Heat kernels and Dirac operators,
  Springer, Berlin-Heidel\-berg-New York (1992)

\Quelle\BH[BH]\Journal
  A.~Borel, F.~Hirzebruch:
  {Characteristic Classes and Homogeneous Spaces, II},
  Amer.\ J.\ Math.~81 (1959), 315--382

\Quelle\Bott[Bo]\ImBuch
  R.~Bott:
  The Index Theorem for Homogeneous Differential Operators,
  in S.~S. Cairns: {Differential and Combinatorial Topology,
  a Symposium in Honor of Marston Morse},
  Princeton Univ. Press (1965), 167--186

\Quelle\Bunke[Bu]\Journal
  U.~Bunke:
  The spectrum of the Dirac operator on the hyperbolic space,
  Math.\ Nachr.~153 (1991), 179--190

\Quelle\CG[CG]\Journal
  M.~Cahen, S.~Gutt:
  Spin Structures on Compact Simply Connected
  {R}iemannian Symmetric Spaces,
  Simon Stevin~62 (1988), 209--242

\Quelle\Campo[C]\Journal
  R.~Camporesi:
  The spinor heat kernel in maximally symmetric spaces,
  Comm.\ Math.\ Phys.~148 (1992), 283--308

\Quelle\CH[CH]\Journal
  R.~Camporesi, A.~Higuchi:
  On the eigenfunctions of the Dirac
  operator on spheres and real hyperbolic spaces,
  J.\ Geom.\ Phys.~20 (1996), 1--18

\Quelle\CM[CM]\Journal
  A.~Connes, H.~Moscovici:
  The $L^2$-Index Theorem for Noncompact Homogeneous Spaces of Lie Groups,
  Ann.\ of Math.~115 (1982), 291--330

\Quelle\FJO[FO]\Journal
  M.~Flensted-Jensen, K.~Okamoto:
  An explicit construction of the $K$-finite vectors in the discrete series
  for an isotropic semisimple symmetric space,
  Mem.\ Soc.\ Math.\ France, Nouv.\ s\'er.~15 (1984), 157--199

\Quelle\GV[GV]\Journal
  E.~Galina, J.~Vargas:
  {Eigenvalues and eigenspaces for the twisted Dirac operator over\/
  $\SU(N,1)$ and \penalty-100$\Spin(2N,1)$},
  Trans.\ Amer.\ Math.\ Soc.~345 (1994), 97--113

\Quelle\Hel[H]\Buch
  S.~Helgason:
  {Differential {G}eometry, {L}ie {G}roups and {S}ymmetric {S}paces},
  Academic Press, New York (1978)

\Quelle\HS[HS]\Journal
  F.~Hirzebruch, P.~Slodowy:
  {Elliptic genera, involutions and homogeneous
  spin manifolds},
  Geom.\ Ded.~35 (1990), 309--343

\Quelle\Knapp[K]\Buch
A.~W.~Knapp:
Representation Theory of Semisimple Groups,
Princeton Univ.\ Press, Princeton, N.\ J.\ (1986)

\Quelle\LM[LM]\Buch
H.~B.~Lawson, M.-L.~Michelsohn:
Spin Geometry,
Princeton Univ.\ Press, Princeton, N.\ J.\ (1989)

\Quelle\Parth[P]\Journal
K.~R.~Parthasarathy:
Dirac operator and the discrete series,
Ann.\ of Math.~96 (1972), 1--30

\Quelle\Seifarth[S]\Preprint
S.~Seifarth:
The discrete spectrum of the Dirac operators on certain symmetric spaces,
Preprint No.~25 of the IAAS, Berlin \yr 1992

\Quelle\Wolf[W]\Journal
J.~A.~Wolf:
Essential self-adjointness for the Dirac operator and its square,
Indiana Univ.\ Math.\ J.~22 (1972/73), 611--640

\endRefs

\enddocument